\documentclass[12pt,oneside,reqno]{amsart}

\usepackage{amssymb}
\usepackage{}
\usepackage{txfonts}
\usepackage{bbm}
\usepackage{cases}
\usepackage{amsmath}
\usepackage{graphicx}
\usepackage{mathrsfs}
\usepackage{stmaryrd}
\usepackage{amsfonts}
\usepackage{hyperref}
\usepackage{enumerate,amsmath,amssymb,amsthm}

\numberwithin{equation}{section}

\newcommand{\be}{\begin{eqnarray}}
\newcommand{\ee}{\end{eqnarray}}
\newcommand{\ce}{\begin{eqnarray*}}
\newcommand{\de}{\end{eqnarray*}}
\newtheorem{theorem}{Theorem}[section]
\newtheorem{lemma}[theorem]{Lemma}
\newtheorem{remark}[theorem]{Remark}
\newtheorem{definition}[theorem]{Definition}
\newtheorem{proposition}[theorem]{Proposition}
\newtheorem{Examples}[theorem]{Example}
\newtheorem{corollary}[theorem]{Corollary}

\def\eps{\varepsilon}

\def\e{\mathrm{e}}

\def\a{\alpha}

\def\u{\mathbf{u}}

\def\p{\partial}
\def\d{\delta}

\def\[{{\Big[}}
\def\]{{\Big]}}
\def\<{{\langle}}
\def\>{{\rangle}}
\def\({{\Big(}}
\def\){{\Big)}}

\def\bx{{\mathbf{x}}}

\def\dif{{\mathord{{\rm d}}}}

\def\min{{\mathord{{\rm min}}}}

\def\no{\nonumber}
\def\={&\!\!=\!\!&}

\def\cJ{{\mathcal J}}

\def\cL{{\mathcal L}}

\def\sS(\R^d){{\mathcal S}}
\def\cT{{\mathcal T}}

\def\mE{{\mathbb E}}

\def\mI{{\mathbb I}}

\def\mL{{\mathbb L}}

\def\mN{{\mathbb N}}

\def\mP{{\mathbb P}}
\def\mQ{{\mathbb Q}}
\def\mR{{\mathbb R}}
\def\mS{{\mathbb S}}

\def\1{{\mathbf{1}}}

\def\sB{{\mathscr B}}

\def\sF{{\mathscr F}}

\def\sL{{\mathscr L}}

\def\sS{{\mathscr S}}

\def\E{\mathbb E}

\def\geq{\geqslant}
\def\leq{\leqslant}

\def\div{\mathord{{\rm div}}}

\def\eps{\varepsilon}

\def\e{\mathrm{e}}

\def\a{\alpha}

\def\u{\mathbf{u}}

\def\p{\partial}
\def\d{\delta}

\def\[{{\Big[}}
\def\]{{\Big]}}
\def\<{{\langle}}
\def\>{{\rangle}}
\def\({{\Big(}}
\def\){{\Big)}}

\def\bx{{\mathbf{x}}}

\def\dif{{\mathord{{\rm d}}}}

\def\min{{\mathord{{\rm min}}}}

\def\no{\nonumber}
\def\={&\!\!=\!\!&}
\def\bt{\begin{theorem}}
\def\et{\end{theorem}}
\def\bl{\begin{lemma}}
\def\el{\end{lemma}}
\def\br{\begin{remark}}
\def\er{\end{remark}}
\def\bx{\begin{Examples}}
\def\ex{\end{Examples}}
\def\bd{\begin{definition}}
\def\ed{\end{definition}}
\def\bp{\begin{proposition}}
\def\ep{\end{proposition}}
\def\bc{\begin{corollary}}
\def\ec{\end{corollary}}

\def\geq{\geqslant}
\def\leq{\leqslant}

\def\div{\mathord{{\rm div}}}

\def\si{\sigma}
 \def\R{\mathbb R}
 \def\R{\mathbb R}

\def\<{\langle} \def\>{\rangle}  
    
\def\d{\text{\rm{d}}}   
  \def\si{\sigma} 
 \def\beq{\begin{equation}}  
 
\def\e{\text{\rm{e}}}

\allowdisplaybreaks

\begin{document}
\title{Well-posedness of supercritical SDE driven by L\'evy processes with irregular drifts}
\author{Zhen-Qing Chen,   Xicheng Zhang and Guohuan Zhao}

\address{Zhen-Qing Chen:
Department of Mathematics, University of Washington, Seattle, WA 98195, USA\\
Email: zqchen@uw.edu
 }

\address{Xicheng Zhang:
School of Mathematics and Statistics, Wuhan University,
Wuhan, Hubei 430072, P.R.China\\
Email: XichengZhang@gmail.com
 }

\address{Guohuan Zhao:
Applied Mathematics, Chinese Academy of Science,
Beijing, 100081, P.R.China\\
Email: zhaoguohuan@gmail.com
 }

\thanks{
Research of ZC is partially supported by  Simons Foundation grant 52054 and NNSFC grant  11731009. 
 Research of XZ is partially supported by NNSFC grant of China (No. 11731009). 
Research of GH is partially supported by National Postdoctoral Program for Innovative Talents (201600182) of China.}

\begin{abstract}
In this paper, we study the following time-dependent stochastic differential equation (SDE) in $\mR^d$:  
$$\dif X_{t}= \si_t(X_{t-}) \d Z_t + b_t(X_{t})\d t,
\quad X_{0}=x\in\R^d,
$$
where $Z$ is a $d$-dimensioanl nondegenerate $\alpha$-stable-like process with $\a\in(0,2)$ (including cylindrical case), 
and uniform in $t\geq 0$,   
$x\mapsto \si_t(x): \mR^d\to\mR^d\otimes\mR^d$ is Lipchitz and uniformly elliptic 
and $x\mapsto b_t (x)$ is $\beta$-order H\"older continuous with $\beta\in(1-\alpha/2,1)$.
Under these  assumptions,
we show the above SDE has a unique strong solution for every starting point $x \in \R^d$. 
When $\sigma_t (x)=\mI_{d\times d}$, the $d\times d$ identity matrix, 
our result in particular gives  an affirmative answer to the open problem in \cite{Pr2}.


\no indent 
 {\bf AMS 2010 Mathematics Subject Classification:}  Primary: 60H10, 35R09; Secondary: 60G51
\end{abstract}

\maketitle

\section{Introduction}

The main purpose of this paper is to establish the strong well-posedness for  a class of supercritical stochastic differential equations
driven by non-degenerate $\alpha$-stable processes, and with time-dependent H\"older drift $b$. More precisely, we are mainly 
concerned 
 with the following SDE:
\begin{equation}\label{SDE}
\d X_{t}=\sigma_t(X_{t-})\d Z_t + b_t(X_{t})\d t,
\quad X_{0}=x\in\R^d,
\end{equation}
where $\sigma:\mR_+\times\mR^d\to\mR^d\otimes\mR^d$ and $b:\mR_+\times\mR^d\to\mR^d$ 
are two Borel measurable functions, and 
$Z$ is a pure jump L\'evy  process with L\'evy measure $\nu$ whose  characteristic function $\phi(x)$  
is given by 
$$
\phi (\xi )=-\int_{{\mR^d}\setminus \{0\}} \left( \mathrm{e}^{i\xi \cdot z} -1 - i \xi  \cdot z 1_{|z |\leq 1}\right) \nu (\dif  z).
$$
When $\sigma_t (x)$ and $b_t (x)$ 
are Lipschitz continuous in $x\in \mR^d$, 
it is well known
that by applying Picard's iteration method and by first removing large jumps of $Z$, one can show
that SDE \eqref{SDE} has a unique strong solution. This paper is concerned with the strong existence and
strong uniqueness solution to SDE \eqref{SDE} when $b_t (x)$ is not Lipschitz continuous in $x$. 

To state our 
condition on L\'evy measure $\nu$, for $\alpha\in(0,2)$, 
 denote by $\mL^{(\alpha)}_{non}$ the space 
of all non-degenerate $\alpha$-stable measures $\nu^{(\alpha)}$, that is,
\begin{align}
\nu^{(\alpha)}(A)=\int^\infty_0\left(\int_{\mS^{d-1}}\frac{1_A (r\theta)\Sigma(\dif\theta)}{r^{1+\alpha}}\right)\dif r,\quad A\in\sB(\mR^d),\label{Eq4}
\end{align}
where $\Sigma$ is a finite measure over the  unit sphere $\mS^{d-1}$ in $\mR^d$ with
 \begin{align} \label{Spe1}
\int_{\mS^{d-1}}|\theta_0\cdot\theta|^\alpha\Sigma(\dif\theta) >0 \quad \hbox{for every } \theta_0\in\mS^{d-1}.  
\end{align} 
Since the left hand side of  the above is a continuous function in $\theta_0 \in \mS^{d-1}$,  
condition \eqref{Spe1} is equivalent to 
$$
\inf_{ \theta_0\in\mS^{d-1}}\int_{\mS^{d-1}}|\theta_0\cdot\theta|^\alpha\Sigma(\dif\theta)>0.  
$$
We assume that there are $\nu_1,\nu_2\in\mL^{(\alpha)}_{non}$,  $\beta \in (1-\alpha/2, 1)$ and $\Lambda>0$ so that
for all $t\geq 0$ and $x, y\in \mR^d$, 
\begin{equation}\label{Con0}
\nu_1(A)\leq \nu(A)\leq\nu_2(A) \quad \hbox{for }   A\in\sB(B_1(0)), 
\end{equation}
\begin{equation}\label{Con1}
|b_t(x)|\leq \Lambda \quad \hbox{and} \quad |b_t(x)-b_t(y)|\leq \Lambda |x-y|^\beta,
\end{equation}
\begin{equation}\label{Con2}
\Lambda^{-1} |\xi|\leq |\si_t(x) \xi|\leq \Lambda |\xi|,
\quad 
\|\sigma_t(x)-\sigma_t(y)\|\leq \Lambda |x-y|,
\end{equation}
where $\|\cdot\|$ denotes the Hilbert-Schmidt norm of a matrix, and $|\cdot|$ denotes the Euclidean norm.
We call a pure jump   L\'evy process $Z$  whose L\'evy measure $\nu$  satisfies condition \eqref{Con0}  
an  $\alpha$-stable-like L\'evy process.
The following is the main result of this paper.

\bt\label{Th11}
Under conditions \eqref{Con0}, \eqref{Con1} and \eqref{Con2}, for each $x\in\mR^d$, there is a unique strong solution to SDE \eqref{SDE}. 
\et

\br
If $Z_t=(Z^{(1)}_t,\cdots, Z^{(d)}_t)$ is a cylindrical $\alpha$-stable process, that is, each component is an independent copy of a 
non-degenerate one-dimensional (possibly asymmetric) 
$\alpha$-stable process, 
then condition \eqref{Spe1} is satisfied. Note that in this case,
 the L\'evy measure of $Z_t$ is singular with respect to the Lebesgue measure.
\er

\br\label{Rek1.3} 
By a standard localization method, if   \eqref{Con1} and \eqref{Con2} are assumed to hold 
on each ball $B(0, R):=\{x\in \mR^d: |x| <R\}$ with $\Lambda$ depends on $R$,
then for each $x\in\mR^d$,
there exists a unique strong solution to SDE \eqref{SDE} up to the explosion time $\zeta$
with $\lim_{t\uparrow\zeta}X_t=\infty$.
\er

The study of weak and strong well-posedness of SDE \eqref{SDE} with irregular coefficients has a long history and there is a large amount of 
literatures devoted to this topic especially when $Z$ is a Brownian motion. 
When $Z$ is a standard $d$-dimensional Brownian motion, $\sigma_t=\mI_{d\times d}$  and $b$ is bounded measurable,
Veretennikov \cite{Ve} proved that SDE \eqref{SDE} has a unique strong solution,
which extended a result of Zvonkin \cite{Zv} 
in one-dimension. 
Using  Girsanov's transformation and  results from PDEs, Krylov and R\"ockner \cite{Kr-Ro} 
obtained the existence and uniqueness of strong solutions to SDE \eqref{SDE} when $\sigma$ is the identity matrix and $b$ satisfies
$$
\|b\|_{L^q_T(L^p(\R^d))}:=\left[\int^T_0 \left(\int_{\mR^d}|b_t(x)|^p\dif x\right)^{q/p}\dif t\right]^{1/q}<\infty, 
\quad  \frac{2}{q}+\frac{d}{p}<1.
$$
These results have been extended to SDEs 
with Sobolev diffusion coefficients and singular drifts in \cite{Zh0, Zh2} by using Zvonkin's idea.

\medskip

However, things become quite different when $Z$ is a pure jump L\'evy process. 
For one-dimensional case, 
Tanaka, Tsuchiya and Watanabe \cite{Ta-Ts-Wa} proved that if $Z$ is a symmetric  $\alpha$-stable process with $\alpha\in[1,2)$,
$\sigma_t (x) \equiv 1$ and $b_t(x)=b(x)$ is 
bounded measurable, then pathwise uniqueness holds for SDE \eqref{SDE}. 
They further show that 
 if $\alpha\in(0,1)$, and even if $b$ is H\"older continuous, the pathwise uniqueness may fail. 
 For one-dimensional multiplicative noise case where $\sigma_t (x)=\sigma (x)$, see \cite{BBC} and \cite[Theorem 1]{Ko}.
For multidimensional case, Priola \cite{Pr1} proved pathwise uniqueness for \eqref{SDE} when  
$\sigma_t (x) =\mI_{d\times d}$, 
$Z$ is a non-degenerate symmetric but possibly non-isotropic $\alpha$-stable process with $\alpha\in[1,2)$ 
and $ b_t (x)=b(x) 
\in C^\beta(\R^d)$ with $\beta\in (1- \alpha/{2},1)$  is time-independent. 
Note that in this case, the infinitesimal generator corresponding to the
solution $X$ of \eqref{SDE} is
$\sL^{(\alpha)}+b\cdot \nabla$. Here $\sL^{(\alpha)}$ is the infinitesimal generator
of the L\'evy process ${Z}$, which is a nonlocal operator of order $\alpha$.
When $\alpha>1$, $\sL^{(\alpha)}$
is the dominant term, which is called the subcritical case.
When $\alpha \in (0, 1)$, the  gradient $\nabla$ is of higher order than
the nonlocal operator $\sL^{(\alpha)}$
so the corresponding SDE \eqref{SDE} is called  supercritical. The critical case corresponds to $\alpha=1$.
Priola's result was extended to drift $b$ in some
fractional Sobolev spaces in the subcritical case in Zhang \cite{Zh1}
and to more general L\'evy processes
in the subcritical and critical cases in Priola \cite{Pr2}. 
Recently, for a large class of L\'evy processes, Chen, Song and Zhang in \cite{Ch-So-Zh} established strong existence and pathwise uniqueness for SDE \eqref{SDE}  when $\sigma_t (x) =\mI_{d\times d}$  and $b_t (x)  $ 
is time-dependent, H\"older continuous in $x$. Therein,
the authors not only extend the main results of \cite{Pr1} and \cite{Pr2} for the subcritical and critical case ($\alpha \in [1,2)$) to more general L\'evy processes and time-dependent drifts $b_t\in L^\infty_T (C^\beta(\R^d))$  
with $\beta\in(1-\frac{\a}{2},1)$, 
but also establish strong existence and pathwise uniqueness for the supercritical case  ($\alpha\in(0,1)$) with $b\in L^\infty_T(C^{\beta}(\R^d))\ (\beta\in (1-\frac{\a}{2},1))$. 
It partially answers an open question posted in \cite{Pr1} on the pathwise well-posedness of SDE \eqref{SDE} in the supercritical case.  
 However, when $Z$ is a  {\em cylindrical} $\a-$stable process $Z$, 
 the result of \cite{Ch-So-Zh} requires $\a>2/3$. 
As mentioned in \cite{Ch-So-Zh}, it is a quite interesting question  whether the constraint $\a> 2/3$ can be dropped. 
 Theorem \ref{Th11} of this paper not only gives  an affirmative answer to the above question but moreover it is
 done for the multiplicative noise setting.  
We remark that except in the one-dimensional case,  almost all the known results  in literature on strong well-posedness of
SDE \eqref{SDE} driven by pure jump L\'evy process $Z$ requires $\sigma_t (x)=\mI_{d\times d}$.
 
\medskip

We now describe the approach of this paper. 
As usual, to study the strong well-posedness of SDE \eqref{SDE}, 
we shall use Zvonkin's transform, which requires a deep understanding for the following nonlocal PDE (Kolmogorov's equation): 
\begin{align}\label{PDE} 
\partial_tu =\sL_t u+b\cdot\nabla u-\lambda u+f \hbox{ with } u(0, x)= 0, 
\end{align}
where
$$
\sL_t u(x):=\int_{\mR^d} [u(x+\sigma_t(x)z)-u(x)-1_{\{ |z|\leq 1\}} \sigma_t(x)z\cdot\nabla u(x)]\nu(\dif z).
$$
We mention that when $\sL_t$ is the usual fractional Laplacian $\Delta^{\alpha/2}:=-(-\Delta)^{\alpha/2}$ with $\alpha\in(0,1)$, 
that is, when 
$\nu(\dif z)=|z|^{-d-\alpha}\dif z$ and $\sigma=\mI_{d\times d}$ in the above definition,
and $b\in L^\infty([0,T];C^\beta)$ with $\beta\in(1-\alpha,1)$, Silvestre \cite{Si2} obtained the following a priori interior 
estimate 
for  any solution $u$ of \eqref{PDE}: 
$$
\|u\|_{L^\infty([0,1];C^{\alpha+\beta}(B_1))} 􏰂\leq C􏰉\Big(\|u\|_{L^\infty([0,2]\times B_2) }+ \|f\|_{L^\infty([0,2];C^\beta(B_2))}\Big),
$$
where $B_r:=\{x\in\mR^d: |x|<r\}$.
Such an interior estimate suggests that one could solve the supercritical SDE \eqref{SDE} uniquely when $Z$ is a rotationally symmetric $\alpha$-stable process with $\alpha\in(0,1)$ and $b\in L^\infty([0,T];C^\beta)$
with $\beta\in(1-\alpha/2, 1)$ (see \cite{Pr2}). 
However the approach of \cite{Si2} strongly depends on realizing the fractional Laplacian in $\R^d$ 
as the boundary trace of an elliptic operator in upper half space of $\R^{d+1}$. 
Extending Silvestre's argument to general  $\alpha$-stable-type operators would be very hard, if 
 not impossible at all. So new ideas are needed for the study of SDE \eqref{SDE} with general L\'evy process $Z$ 
 and  variable diffusion matrix $\sigma_t(x)$.

\medskip 

Our approach of studying \eqref{PDE} is based on the Littlewood-Paley decomposition and some Bernstein's type inequalities. 
This approach allows us to handle a large class of L\'evy's type operator in a 
unified way, including  L\'evy's type operators
with singular L\'evy measures, see Theorem \ref{Priori} below. 
When $\sigma_t(x)=\sigma_t$ is spatially independent and 
real part of the symbol $\psi_t(\xi)$ of $\sL_t$  
(that is,  $\sF (\sL_tf )(\xi)= \psi_t(\xi ) \sF(f)(\xi)$, where $\sF (f)$ denotes the Fourier transform of $f$)
 is bounded from above by 
 $-c_0|\xi|^\alpha$, we show the following a priori estimate for \eqref{PDE}:
 for every $p>d/(\alpha+\beta-1)$, there is a constant $C>0$ depending only on $T,d,p,\alpha,\beta$ and $\|b\|_{L^\infty([0,T]; B^\beta_{p,\infty})}$
 so that 
$$
\|u\|_{L^\infty([0,T];B^{\alpha+\beta}_{p,\infty})} 􏰂\leq 􏰉C \|f\|_{L^\infty([0,T];B^\beta_{p,\infty})},
$$ 
where $B^{\beta}_{p,\infty}$ is the usual Besov space (see Definition \ref{Def2} below). 
The above  a priori estimate is the key in our solution to  the pathwise well-posedness problem of SDE \eqref{SDE}
when $\sigma_t (x)$ is spatially independent. 
The  general case with variable  coefficient $\sigma_t(x)$ is much more delicate. 
First of all, in general 
$$
x\mapsto \int_{|z|>1}f(x+\sigma_t(x) z)\nu(\dif z)
\mbox{ may not be smooth even if $f (x)$ and $\sigma_t(x)$ are smooth}.
$$
Thus  to treat the general case, we have to  first remove the large jumps. 
Next we need to
impose a small condition on the oscillation 
of $\sigma$ by using a perturbation argument, see Theorem \ref{Th34} 
below, which will be removed later through a localization and patching together procedure. 
\medskip 

This paper is organized as follows: In Section 2, we recall some well-known facts from Littlewood-Paley theory, 
 in particular, the Bony decomposition and  Bernstein's inequalities,   and establish a useful commutator estimate. 
In Section 3, we study the nonlocal advection equation \eqref{PDE} with irregular drift $b$,
and obtain some a priori estimates in Besov spaces. In Section 4, we prove our main theorem by Zvonkin's transform and 
a suitable patching together technique.

\medskip

We close this section by mentioning 
some conventions used throughout this paper:
We use $:=$ as a way of definition. For $a, b\in \mR$, $a\vee b:= \max \{a, b\}$ and $a\wedge b:=\min \{a, b\}$,
and on $\mR^d$, $\nabla:=(\frac{\partial}{\partial x_1}, \dots, \frac{\partial}{\partial x_d})$ and 
$\Delta:= \sum_{k=1}^d \frac{\partial^2}{\partial x_k^2}$.
The letter $c$ or $C$ with or without subscripts stands for an unimportant constant, whose value may change in difference places.
We use $A\asymp B$ to denote that $A$ and $B$ are comparable up to a constant, and use $A\lesssim B$ to denote $A\leq C B$ 
for some constant $C$.

\section{Preliminary}

 In this section, we recall  some basic facts from Littlewood-Paley theory, especially Bernstein's inequalities (see \cite{Ba-Ch-Da}). 
We then establish a commutator estimate, which plays an important role in our approach.  

Let $\sS(\R^d)$ be the Schwartz space of all rapidly decreasing functions, and $\sS'(\R^d)$ the dual space of $\sS(\R^d)$ 
called Schwartz generalized function (or tempered distribution) space. Given $f\in\sS(\R^d)$,
  let $\sF f=\hat f$ be the Fourier transform of $f$ defined by
$$
\hat f(\xi):=(2\pi)^{-d/2}\int_{\mR^d}\e^{-\mathrm{i}\xi\cdot x} f(x)\dif x.
$$
For $R, R_1, R_2\geq0$ with $R_1<R_2$, we shall denote
$$
B_R:=\{x\in\mR^d: |x|\leq R\},\ \ D_{R_1, R_2}:=\{x\in\mR^d: R_1\leq x\leq R_2\}.
$$
The following simple fact will be used frequently: Let $f, g\in\sS'(\R^d)$ be two tempered distributions with supports in $B_{R_0}$ and
$D_{R_1,R_2}$ respectively. Then 
\begin{align}\label{EE6}
\mbox{supp$f*g\subset D_{(R_1-R_0)\vee 0, R_2+R_0}$.}
\end{align}
Let $\chi:\mR^d\to[0,1]$ be a smooth radial function with 
$$
\chi(\xi)=1,\ |\xi|\leq 1,\ \chi(\xi)=0,\ |\xi|\geq 3/2.
$$
Define
$$
\varphi(\xi):=\chi(\xi)-\chi(2\xi).
$$
It is easy to see that $\varphi\geq 0$ and supp $\varphi\subset B_{3/2}\setminus B_{1/2}$ and
\begin{align}\label{EE1}
\chi(2\xi)+\sum_{j=0}^k\varphi(2^{-j}\xi)=\chi(2^{-k}\xi)\stackrel{k\to\infty}{\to} 1.
\end{align}
In particular, if $|j-k|\geq 2$, then
$$
\mathrm{supp} \left[ \varphi(2^{-j}\cdot) \right] \cap\mathrm{supp} \left[ \varphi(2^{-k}\cdot) \right] =\emptyset.
$$
From now on we shall fix such $\chi$ and $\varphi$, and introduce the following definitions. 
\bd\label{Def2}
The dyadic block operator 
$\Pi_j$ is defined by
$$
\Pi_j f:=
\left\{
\begin{array}{ll}
\sF^{-1}(\chi(2\cdot) \sF f), & j=-1, \\
\sF^{-1}(\varphi(2^{-j}\cdot) \sF f),& j\geq 0.
\end{array}
\right.
$$
For $s\in\mR$ and $p,q\in[1,\infty]$, the Besov space $B^s_{p,q}$ is defined as the set of all $f\in\sS'(\R^d)$ with
$$
\|f\|_{B^s_{p,q}}:=1_{\{q<\infty\}}\left(\sum_{j\geq -1}2^{jsq}\|\Pi_j f\|_p^q\right)^{1/q}+1_{\{q=\infty\}}\left(\sup_{j\geq -1}2^{js}\|\Pi_j f\|_p\right)<\infty,
$$
where $\|\cdot\|_p$ denotes the usual $L^p$-norm in $\mR^d$. 
\ed

Some literature, e.g., \cite{Ba-Ch-Da, Ch-Mi-Zh},   uses notation $\Delta_j$  for the dyadic block operator $\Pi_j$ defined above. 
We choose to use notation $\Pi_j$ in this paper out of two considerations: 
the dyadic block operator is a projection operator in the $L^2$-space and we want to avoid possible confusion
with Laplacian operator $\Delta$.
 
Let $H^s_p:=(I-\Delta)^{-s/2}(L^p)$ be the usual Bessel potential space with norm
 $$
\|f\|_{H^s_p}:=\|(I-\Delta)^{s/2} f\|_p. 
$$
Note that  if $s\geq 0$, $\|f\|_{H^s_p}\asymp \|f\|_p+\|(-\Delta)^{s/2}f\|_p$. 
It should be observed that if $s>0$ is not an integer, then Besov space $B^s_{\infty,\infty}$ is just the usual H\"older space $C^s$.
Moreover, Besov spaces have the following embedding relations: For any $s,s',s''\in\mR$ and $p,p',q,q'\in[1,\infty]$ with 
$$
p\leq p',  \ q\leq q',\  s<s'' \ \hbox{ with } s-d/p=s'-d/p',  
$$ 
it holds that (cf. \cite{Be-Lo})
\begin{align}\label{Emm}
B^{s''}_{p,1}\subset H^{s''}_p\subset B^{s''}_{p,\infty}\subset B^s_{p,q}\subset B^{s'}_{p',q'}.
\end{align}
Let $h=\sF^{-1} \chi$ be the inverse Fourier transform of $\chi$. Define
$$
h_{-1}(x):=\sF^{-1} \chi(2\cdot)(x)=2^{-d}h(2^{-1}x)\in\sS(\mR^d),
$$
and for $j\geq 0$,
\begin{align}\label{EE7}
h_j(x):=\sF^{-1}\varphi(2^{-j}\cdot)(x)=2^{jd}h(2^jx)-2^{(j-1)d}h(2^{j-1}x)\in \sS(\mR^d).
\end{align}
By definition it is easy to see that 
\begin{align}\label{EE2}
\Pi_j f(x)=(h_j*f)(x)=\int_{\mR^d}h_j(x-y)f(y)\dif y,\ \ j\geq -1.
\end{align}
The cut-off low frequency operator $S_k$ is defined by
$$
S_kf:=\sum_{j=-1}^{k-1}\Pi_j f=2^{(k-1)d}\int_{\mR^d}h(2^{k-1}(x-y))f(y)\dif y.
$$
It is easy to see that 
\begin{align}\label{ED1}
\|S_kf\|_p\leq\|h\|_1\|f\|_p,\ \ \|S_kf\|_{B^s_{p,q}}\leq \|h\|_1\|f\|_{B^s_{p,q}}.
\end{align}
Moreover, by \eqref{EE1}, one has
\begin{align}
\widehat{S_k f}=\chi(2^{1-k}\cdot)\hat f,\ \ f=\lim_{k\to\infty}S_k f=\sum_{j\geq -1}\Pi_j f.\label{ED4}
\end{align}
For $f,g\in B^s_{p, q}$, define
$$
T_fg=\sum_k S_{k-1}f\Pi_k g,\ \ R(f,g):=\sum_k\sum_{|i|\leq1}\Pi_k f\Pi_{k-i}g.
$$
The following identity 
$$
fg=T_fg+T_g f+R(f,g)
$$
is called the  Bony decomposition of $fg$.

We first recall the following Bernstein's type inequality (cf. \cite{Ba-Ch-Da} and \cite{Ch-Mi-Zh}).

\bl\label{Le13}
(Bernstein's type inequality) 
Let $1\leq p\leq q\leq \infty$. For any $k=0,1,\cdots$ and $\beta\in(-1,2)$, there is a constant $C$ such that for all $f\in\sS'(\mR^d)$ and $j\geq -1$, 
\begin{align}\label{EE9}
\|\nabla^k\Pi_j f\|_q\leq C 2^{(k+d(\frac{1}{p}-\frac{1}{q}))j}\|\Pi_jf\|_p,
\end{align} 
and for any $j\geq 0$, 
\begin{align}\label{EE99}
\|(-\Delta)^{\beta/2}\Pi_j f\|_q\leq C 2^{(\beta+d(\frac{1}{p}-\frac{1}{q}))j}\|\Pi_jf\|_p,
\end{align}
and for any $2\leq p<\infty$, $j\geq 0$ and $\a\in(0,2)$, there is a constant $c>0$ such that for all $f\in\sS'(\mR^d)$, 
\begin{align}\label{New}
\int_{\mR^d}\Big|(-\Delta)^{\alpha/4}|\Pi_j f|^{p/2}\Big|^2\dif x\geq c 2^{\alpha j}\|\Pi_j f\|_p^p. \end{align}
\el

The following commutator estimate plays an important role in this paper.
\bl\label{Le12}
Let $p, p_1, p_2, q_1, q_2\in[1,\infty]$ with $\frac{1}{p}=\frac{1}{p_1}+\frac{1}{p_2}$ and $\frac{1}{q_1}+\frac{1}{q_2}=1$.
For any $\beta_1\in(0,1)$ and $\beta_2\in[-\beta_1,0]$, there is a constant $C>0$ depending only on $d,p,p_1,p_2,\beta_1,\beta_2$ such that 
 $$
\|[\Pi_j, f]g\|_p\leq C2^{-j(\beta_1+\beta_2)}
\left\{
\begin{array}{ll}
\|f\|_{B^{\beta_1}_{p_1,\infty}}\|g\|_{p_2},&\mbox{if}\ \  \beta_2=0,\\
\|f\|_{B^{\beta_1}_{p_1,\infty}}\|g\|_{B^{\beta_2}_{p_2,\infty}},& \mbox{if}\ \  \beta_1+\beta_2>0,\\
\|f\|_{B^{\beta_1}_{p_1,q_1}}\|g\|_{B^{\beta_2}_{p_2,q_2}},& \mbox{if}\ \ \beta_1+\beta_2=0,
\end{array}
\right.
$$
where $[\Pi_j, f]g:=\Pi_j (fg)-f\Pi_jg$. 
\el

\begin{proof}
We first consider the case $\beta_2=0$. In this case , by \eqref{EE2}, 
\begin{align*}
[\Pi_j, f]g(x)=\int_{\mR^d}h_j(y)(f(x-y)-f(x))g(x-y)\dif y.
\end{align*}
 For any $p\in[1,\infty]$ and $s\in(0,1)$, by Theorem 2.36 of \cite{Ba-Ch-Da},
\begin{align}\label{ED8}
\|f(\cdot-y)-f(\cdot)\|_p\leq C|y|^{s}\|f\|_{B^{s}_{p,\infty}}.
\end{align}
Using
H\"older's inequality and \eqref{EE7}, we have
\begin{align}\label{EE5}
\begin{split}
\|[\Pi_j, f]g\|_p&\leq \int_{\mR^d}h_j(y)\|f(\cdot-y)-f(\cdot)\|_{p_1}\|g\|_{p_2}\dif y\\
&\lesssim\|f\|_{B^{\beta_1}_{p_1,\infty}}\|g\|_{p_2}\int_{\mR^d}|h_j(y)|\,|y|^{\beta_1}\dif y\\
&=\|f\|_{B^{\beta_1}_{p_1,\infty}}\|g\|_{p_2}2^{-j\beta_1} \int_{\mR^d}|2h(2y)-h(y)|\,|y|^{\beta_1}\dif y\\
&\lesssim 2^{-j\beta_1} \|f\|_{B^{\beta_1}_{p_1,\infty}}\|g\|_{p_2}.
\end{split}
\end{align}
Next we consider the case $\beta_2\in[-\beta_1,0)$.
By using Bony's decomposition, we can write
\begin{align*}
[\Pi_j, f]g
=[\Pi_j,T_{f}] g+\Pi_j(T_g f)-T_{\Pi_jg} f+\Pi_jR(f,g)-R(f,\Pi_jg).
\end{align*}
Notice that by \eqref{ED1} and \eqref{EE6},
$$
\sF(\Pi_j (S_{k-1}f\Pi_k g))=\varphi_j\cdot (\chi(2^{1-k}\cdot) \hat f)*(\varphi_k\hat g)=0 \ \mbox{ for }\ |k-j|>4.
$$ 
Therefore, 
by \eqref{ED1} and \eqref{EE5} we have
\begin{align*}
\|[\Pi_j,T_{f}] g\|_p&=\Bigg\|\sum_{|k-j|\leq 4}\Big(\Pi_j (S_{k-1}f\Pi_k g)-S_{k-1}f\Pi_j \Pi_k g\Big)\Bigg\|_p\\
&\leq\sum_{|k-j|\leq 4}\big\|[\Pi_j, S_{k-1}f]\Pi_k g\big\|_p\\
&\lesssim 2^{-j\beta_1}\sum_{|k-j|\leq 4}\|S_{k-1}f\|_{B^{\beta_1}_{p_1,\infty}}\|\Pi_k g\|_{p_2}\\
&\lesssim 2^{-j\beta_1}\|f\|_{B^{\beta_1}_{p_1,\infty}}\sum_{|k-j|\leq 4}2^{-k\beta_2}\|g\|_{B^{\beta_2}_{p_2,\infty}}\\
&\lesssim 2^{-j(\beta_1+\beta_2)}\|f\|_{B^{\beta_1}_{p_1,\infty}}\|g\|_{B^{\beta_2}_{p_2,\infty}}.
\end{align*}
Similarly, by H\"older's inequality and $\beta_2<0$, we have
\begin{align*}
\|\Pi_j(T_g f)\|_p&=\Bigg\|\sum_{|k-j|\leq 4}\Pi_j(S_{k-1} g\Pi_k f)\Bigg\|_p\leq \sum_{|k-j|\leq 4}\|\Pi_j(S_{k-1} g\Pi_k f)\|_p\\
&\leq \sum_{|k-j|\leq 4}\|S_{k-1} g\Pi_k f\|_p\leq C\sum_{|k-j|\leq 4}\sum_{m\leq k-2}\|\Pi_m g\,\Pi_k f\|_p\\
&\lesssim \|g\|_{B^{\beta_2}_{p_2,\infty}}\|f\|_{B^{\beta_1}_{p_1,\infty}}\sum_{|k-j|\leq 4}\sum_{m\leq k-2}2^{-m\beta_2}2^{-k\beta_1}\\
&\lesssim \|g\|_{B^{\beta_2}_{p_2,\infty}}\|f\|_{B^{\beta_1}_{p_1,\infty}}2^{-j(\beta_2+\beta_1)},
\end{align*}
and
\begin{align*}
\|T_{\Pi_jg} f\|_p&\leq\sum_{k\geq j-2}\|S_{k-1}\Pi_jg\Pi_k f\|_p\leq\sum_{k\geq j-2}\|\Pi_k f\|_{p_1}\|S_{k-1}\Pi_jg\|_{p_2}\\
&\leq\sum_{k\geq j-2}2^{-k\beta_1}\|f\|_{B^{\beta_1}_{p_1,\infty}}\|\Pi_jg\|_{p_2}
\leq C2^{-j(\beta_1+\beta_2)}\|f\|_{B^{\beta_1}_{p_1,\infty}}\|g\|_{B^{\beta_2}_{p_2,\infty}}.
\end{align*}
Finally, we have
\begin{align*}
\|\Pi_jR(f,g)\|_p&=\Bigg\|\sum_{|i|\leq 1, k\geq j-4}\Pi_j(\Pi_k f\Pi_{k-i} g)\Bigg\|_p
\lesssim \sum_{|i|\leq 1, k\geq j-4}\|\Pi_k f\|_{p_1}\|\Pi_{k-i} g\|_{p_2}\\
&\lesssim \sum_{|i|\leq 1, k\geq j-4}2^{-k(\beta_1+\beta_2)}\Big(2^{k\beta_1}\|\Pi_kf\|_{p_1}\Big)\Big(2^{k\beta_2}\|\Pi_{k-i}g\|_{p_2}\Big)\\
&\lesssim 2^{-j(\beta_1+\beta_2)}
\left\{
\begin{array}{lc}
\|f\|_{B^{\beta_1}_{p_1,\infty}}\|g\|_{B^{\beta_2}_{p_2,\infty}},&\beta_1+\beta_2>0,\\
\|f\|_{B^{\beta_1}_{p_1,q_1}}\|g\|_{B^{\beta_2}_{p_2,q_2}} ,&\beta_1+\beta_2=0,
\end{array}
\right.
\end{align*}
where $\frac{1}{q_1}+\frac{1}{q_2}=1$, and
\begin{align*}
&\|R(f,\Pi_jg)\|_p=\Bigg\|\sum_{|i|\leq 1, |k-j|\leq 1}\Pi_{k-i} f\Pi_{k}\Pi_j g\Bigg\|_p
\lesssim \|f\|_{B^{\beta_1}_{p_1,\infty}}\|g\|_{B^{\beta_2}_{p_2,\infty}}2^{-j(\beta_1+\beta_2)}.
\end{align*}
Combining the above calculations and using $\|f\|_{B^s_{p,\infty}}\leq\|f\|_{B^s_{p,q}}$, we complete the proof. 
\end{proof}

\section{Nonlocal parabolic equations}

In this section we study the solvability and regularity of nonlocal parabolic equation \eqref{PDE} with H\"older drift.
Let $\sigma$ be a constant $d\times d$-matrix
and $\nu$ a   measure on $\mR^d$ such that  
$$
\int_{\mR^d\setminus \{0\}}(|z|^2\wedge 1)\nu(\dif z)<\infty.
$$
We define  a L\'evy-type operator $ \cL^\nu_\sigma$ by
$$
\cL^\nu_\sigma f(x):=\int_{\mR^d}\Big(f(x+\sigma z)-f(x)-1_{\{ |z|\leq 1\}}\sigma z\cdot\nabla f(x)\Big)\nu(\dif z),\ f\in\sS(\mR^d).
$$
By Fourier's transform, we have
$$
\widehat{\cL^\nu_\sigma f}(\xi)=\psi^\nu_\sigma(\xi)\hat f(\xi),
$$
where the symbol $\psi^\nu_\sigma(\xi)$ 
is given by 
$$
\psi^\nu_\sigma(\xi)=\int_{\mR^d} \left(\e^{\mathrm{i}\xi\cdot \sigma z}-1-1_{\{ |z|\leq 1\}} \mathrm{i}\sigma z\cdot\xi \right)\nu(\dif z).
$$
Now let $\sigma_t(x): \mR_+\times \mR^d\to\mR^d\otimes\mR^d$ be a Borel measurable function. 
Define a time-dependent L\'evy-type operator 
$$
\sL_tf(x):=\cL^\nu_{\sigma_t(x)}f(x).
$$
In this section, for  $\lambda \geq 0$,  we study
the solvability of the following equation with Besov drift $b_t(x):\mR_+\times\mR^d\to\mR^d$,
\begin{align}\label{EE8}
\p_t u=(\sL_t u-\lambda) u+b\cdot\nabla u+f \hbox{ with }  u(0)=0. 
\end{align}
For a space-time function $f:\mR_+\times\mR^d\to\mR$ and $T>0$, define 
$$
\|f\|_{L^\infty_T(B^s_{p,q})}: =\sup_{t\in[0,T]} \|f(t,\cdot)\|_{B^s_{p,q}}. 
$$

 \subsection{Constant coefficient case: $\sigma_t(x)=\sigma_t$}
 
In this subsection we consider equation \eqref{EE8} 
with time dependent constant coefficient $\sigma_t(x)=\sigma_t$.
First of all, we establish the 
following Bernstein's type inequality for nonlocal operator $\cL^\nu_\sigma$, which plays a crucial role in the sequel.
\bl\label{Le13}
Suppose {\rm Re}$(\psi^\nu_\sigma(\xi))\leq -c_0|\xi|^\alpha$ for some $c_0>0$. Then for any 
$p>2$, there is a constant 
$c_p=c(c_0,p)>0$ such that for $j=0,1,\cdots$,
\begin{align}\label{Bernstein2}
 \int_{\mR^d}|\Pi_j f|^{p-2} ( \Pi_j f) \cL^\nu_\sigma \Pi_j f\dif x\leq -c_p 2^{\alpha j}\|\Pi_j f\|_p^p,
  \end{align}
and for $j=-1$,
$$
 \int_{\mR^d}|\Pi_{-1} f|^{p-2} (\Pi_{-1} f) \cL^\nu_\sigma \Pi_{-1} f\dif x\leq 0.
 $$
\el

\begin{proof}
For $p\geq 2$, by the elementary inequality 
 $|r|^{p/2}-1\geq \frac{p}{2}(r-1)$ for $r\in\mR$, 
we have
$$
 |a|^{p/2}-|b|^{p/2}\geq \tfrac{p}{2}(a-b)b|b|^{p/2-2} ,\  \ a,b\in\mR.
 $$
Letting $g$ be a smooth function, by definition we have
\begin{align*}
\cL^\nu_\sigma |g|^{p/2}(x)&=\int_{\mR^d}\Big(|g(x+\sigma z)|^{p/2}-|g(x)|^{p/2}-1_{|z|\leq 1}\sigma z\cdot\nabla_x|g(x)|^{p/2}\Big)\nu(\dif z)\\
 &\geq\frac{p}{2}|g(x)|^{p/2-2}g(x)\int_{\mR^d}\Big(g(x+\sigma z)-g(x)-1_{|z|\leq 1}\sigma z\cdot\nabla g(x)\Big)\nu(\dif z)\\
 &=\frac{p}{2}|g(x)|^{p/2-2}g(x)\cL^\nu_\sigma g(x).
 \end{align*}
 Multiplying both sides by $|g|^{p/2}$ and then integrating in $x$ over $\mR^d$, by Plancherel's formula, we obtain
\begin{align*}
\int_{\mR^d}|g|^{p-2}g\cL^\nu_\sigma g\dif x&\leq\frac{2}{p}\int_{\mR^d}|g|^{p/2}\cL^\nu_\sigma |g|^{p/2}\dif x
=\frac{2}{p}\int_{\mR^d}|\widehat{|g|^{p/2}}(\xi)|^2\psi^\nu_\sigma(\xi)\dif\xi\\
&=\frac{2}{p}\int_{\mR^d}|\widehat{|g|^{p/2}}(\xi)|^2\mathrm{Re}(\psi^\nu_\sigma(\xi))\dif\xi
\leq -\frac{2c_0}{p}\int_{\mR^d}|\widehat{|g|^{p/2}}(\xi)|^2|\xi|^\alpha\dif\xi\\
&\leq -\frac{2c_0}{p}\int_{\mR^d}|(-\Delta)^{\alpha/4}|g|^{p/2}|^2\dif x,
\end{align*}
which in turn gives the desired estimate by taking $g=\Pi_j f$ and \eqref{New}.
\end{proof}

Now we can state our main result of this subsection.

\bt\label{Priori}
Let $\beta\in(0,1)$ and $\alpha\in(0,2)$ with $\alpha+\beta>1$. Let $T>0$ and 
$p\in(\frac{d}{\alpha+\beta-1}\vee 2 ,\infty)$.
Suppose that for some $c_0>0$ and all $t\in[0,T]$,
\begin{align}\label{ED2}
\mbox{{\rm Re}$(\psi^\nu_{\sigma_t}(\xi))\leq -c_0|\xi|^\alpha$},\ \ \xi\in\mR^d,
\end{align}
and $b=b_1+b_2$ with
$$
b_1\in L^\infty_T(B^\beta_{p,\infty}) \quad \hbox{and} \quad  b_2\in L^\infty_T(B^\beta_{\infty,\infty}).
$$ 
Then for any  $\gamma\in[0,\beta]$ and $f\in L^\infty_T(B^\gamma_{p,\infty})$, there exists a unique 
solution $u\in L^\infty_T (B^{\alpha+\gamma}_{p,\infty})$ to equation \eqref{EE8} in the weak sense, i.e. for all $\varphi\in \sS(\mR^d)$,
$$
\<u(t),\varphi\>=\int^t_0\<u,\ (\sL^*_s -\lambda)\varphi\>\dif s+\int^t_0\<b\cdot\nabla u+f,\varphi\>\dif s,
$$
where $\<u,\varphi\>:=\int_{\mR^d} u\varphi\dif x$ and $\sL^*_s$ is the adjoint operator of $\sL_s$. Moreover, there is a constant 
$C>0$ depending only on $T, d,p,\alpha,\beta,\gamma$ and $\|b_1\|_{L^\infty_T(B^\beta_{p,\infty})}$, $\|b_2\|_{L^\infty_T(B^\beta_{\infty,\infty})}$
 such that for  all $\lambda\geq 0$,
\begin{equation}\label{Hightorder}
\|u\|_{L^\infty_T (B^{\alpha+\gamma}_{p,\infty})}\leq C \|f\|_{L^\infty_T(B^{\gamma}_{p,\infty})},
\end{equation}
and for any $s\in[0,\alpha+\gamma)$, 
\begin{equation}\label{ED7}
\|u\|_{L^\infty_T(B^{s}_{p,\infty})}\leq c_{\lambda}\|f\|_{L_T^\infty (B^{\gamma}_{p,\infty})},
\end{equation}
where $c_{\lambda}=c(\lambda, \a, s, d,p)\rightarrow 0$ as $\lambda\rightarrow \infty$. 
\et

\begin{proof}
(i) We first assume 
$$
b_1, f\in \cap_{s\geq 0}L^\infty_T(B^s_{p,\infty})
\quad \hbox{and} \quad  b_2\in\cap_{s\geq 0}L^\infty_T(B^s_{\infty,\infty}).
$$
Under this assumption, it is well-known 
that  the non-local 
PDE \eqref{EE8} has a unique smooth solution $u$ (see \cite{Zh2}). 
Our main task is to show the a priori estimates \eqref{Hightorder}
and \eqref{ED7}.
Using operator $\Pi_j$ acts on both sides of \eqref{EE8}, we have
\begin{align*}
\p_t \Pi_ju=&(\sL_t-\lambda)\Pi_j u+\Pi_j(b\cdot\nabla u)+\Pi_jf\\
=&(\sL_t-\lambda)\Pi_j u+(b\cdot\nabla \Pi_j u)+[\Pi_j, b\cdot\nabla ]u+\Pi_jf.
\end{align*}
For $p>2$, 
 by the chain rule or multiplying both sides by $|\Pi_j u|^{p-2}\Pi_j u$ and then integrating in $x$, we obtain
\begin{align*}
\frac{\p_t\|\Pi_j u\|^p_p}{p} =& \int_{\mR^d}\Big(|\Pi_j u|^{p-2} (\Pi_j u) \big(\sL_t\Pi_j u+\Pi_j(b\cdot\nabla u)+\Pi_jf-\lambda \Pi_j u\big)\Big)\dif x\\
 =& \int_{\mR^d}|\Pi_j u|^{p-2} (\Pi_j u) \sL_t\Pi_j u\dif x+\int_{\mR^d}|\Pi_j u|^{p-2} (\Pi_j u) \,[\Pi_j,b\cdot\nabla]u\dif x\\
&+\!\!\int_{\mR^d}|\Pi_j u|^{p-2} (\Pi_j u) \,(b\cdot\nabla) \Pi_j u\dif x+\!\!\int_{\mR^d}|\Pi_j u|^{p-2} (\Pi_j u)\Pi_jf\dif x-\lambda\|\Pi_j u\|^p_p\\
 =:&  I^{(1)}_j+I^{(2)}_j+I^{(3)}_j+I^{(4)}_j+I^{(5)}_j.
\end{align*}
For $I^{(1)}_j$, recalling $\sL_t=\cL^\nu_{\sigma_t}$ and by Lemma \ref{Le13}, there is a $c>0$ such that
\begin{align*}
I^{(1)}_{-1}\leq 0,\ \ I^{(1)}_j\leq -c2^{\alpha j}\|\Pi_j u\|^p_p,\ \ j=0,1,2,\cdots.
\end{align*}
For $I^{(2)}_j$, using Lemma \ref{Le12} with
$$
f=b^i,  \quad g=  \p_i u  \quad \hbox{for } i=1,\cdots, d,
$$
and 
$$
\beta_1 = \beta, \  \   \beta_2 = \gamma-\beta, \   \   q_1=\infty \quad \hbox{and} \quad q_2=1, 
$$
by H\"older's inequality and recalling $b=b_1+b_2$, we have for all $j=-1,0,1,\cdots$,
\begin{align*}
I^{(2)}_j&\leq\|[\Pi_j,b\cdot\nabla]u\|_p\|\Pi_j u\|_p^{p-1}
\lesssim 2^{-\gamma j}\Big(\|b_1\|_{B^\beta_{p,\infty}}\|u\|_{B^{1-\beta+\gamma}_{\infty,1}}+\|b_2\|_{B^\beta_{\infty,\infty}}\|u\|_{B^{1-\beta+\gamma}_{p,1}}\Big)
\|\Pi_j u\|_p^{p-1}.
\end{align*}
For $I^{(3)}_j$, note that 
\begin{align*}
I^{(3)}_j&=\int_{\mR^d}((b-S_jb)\cdot\nabla) \Pi_j u\,|\Pi_j u|^{p-2}\Pi_j u\dif x 
  +\int_{\mR^d}(S_jb\cdot\nabla) \Pi_j u\,|\Pi_j u|^{p-2}\Pi_j u\dif x \\
  & =:I^{(31)}_j+I^{(32)}_j.
\end{align*}
By Bernstein's inequality \eqref{EE9}, we have
\begin{align*}
I^{(31)}_j&\leq\sum_{k\geq j}\|(\Pi_kb\cdot\nabla) \Pi_j u\|_p\|\Pi_j u\|^{p-1}_p\\
&\leq\sum_{k\geq j}\Big(\|\Pi_kb_1\|_{p}\|\nabla\Pi_j u\|_{\infty}+\|\Pi_kb_2\|_{\infty}\|\nabla\Pi_j u\|_p\Big)\|\Pi_j u\|^{p-1}_p\\
&\lesssim 2^{(1+d/p)j}\|\Pi_j u\|^{p}_p\sum_{k\geq j}\Big(\|\Pi_kb_1\|_{p}+\|\Pi_k b_2\|_\infty\Big)\\
&\lesssim 2^{(1+d/p-\beta)j}\|\Pi_j u\|^{p}_p \Big(\|b_1\|_{B^\beta_{p,\infty}}+\|b_2\|_{B^\beta_{\infty,\infty}}\Big).
\end{align*}
For $I^{(32)}_j$, we have  by the divergence theorem and \eqref{EE9} again, 
\begin{align*}
I^{(32)}_j&=\frac{1}{p}\int_{\mR^d}(S_jb\cdot\nabla) |\Pi_j u|^p\dif x=-\frac{1}{p}\int_{\mR^d}S_j\div b\, |\Pi_j u|^p\dif x\\
&\leq \frac{1}{p}\|S_j\div b\|_{\infty}\|\Pi_j u\|^p_{p}\leq \frac{1}{p}\sum_{k\leq j}\|\Pi_k \div b\|_{\infty}\|\Pi_j u\|^p_{p}\\
&\lesssim \sum_{k\leq j}2^{k(1+d/p)}\Big(\|\Pi_kb_1\|_{p}+\|\Pi_kb_2\|_{\infty}\Big)\|\Pi_j u\|^p_{p}\\
&\lesssim 2^{(1-\beta+d/p)j}\Big(\|b_1\|_{B^\beta_{p,\infty}}+\|b_2\|_{B^\beta_{\infty,\infty}}\Big)\|\Pi_j u\|^p_{p}.
\end{align*}
Combining the above two estimates, we obtain
\begin{align*}
\p_t\|\Pi_j u\|^p_p/p&\leq-c2^{\alpha j}1_{j\geq 0}\|\Pi_j u\|^p_p-\lambda \|\Pi_j u\|_p^p\\
&\quad+C2^{-\gamma j}\Big(\|u\|_{B^{1-\beta+\gamma}_{\infty,1}}+\|u\|_{B^{1-\beta+\gamma}_{p,1}}\Big)\|\Pi_j u\|_p^{p-1}\\
&\quad+C2^{(1-\beta+d/p)j}\|\Pi_j u\|^p_{p}+C\|\Pi_j u\|^{p-1}_{p}\|\Pi_j f\|_p\\
&\leq-\Big(c2^{\alpha j}1_{j\geq 0}+\lambda -C2^{(1-\beta+d/p)j}\|b\|_{B^\beta_{p,\infty}}\Big)\|\Pi_j u\|^p_p\\
&\quad+C\Big(2^{-\gamma j}\big(\|u\|_{B^{1-\beta+\gamma}_{\infty,1}}
+\|u\|_{B^{1-\beta+\gamma}_{p,1}}\big)+\|\Pi_j f\|_p\Big)\|\Pi_j u\|_p^{p-1}.
\end{align*}
Since $1-\beta+d/p<\alpha$, by dividing both sides by $\|\Pi_j u\|_p^{p-1}$ and using Young's inequality, we get for some $c_0,\lambda_0>0$ and all $j\geq -1$,
$$
\p_t\|\Pi_j u\|_p
\leq-(c_02^{\alpha j}+\lambda-\lambda_0)\|\Pi_j u\|_p
+C2^{-\gamma  j}\big(\|u\|_{B^{1-\beta+\gamma}_{\infty,1}}
+\|u\|_{B^{1-\beta+\gamma}_{p,1}}\big)+C\|\Pi_j f\|_p,
$$
which implies by Gronwall's inequality that for all $j\geq -1$,
\begin{align}
&\|\Pi_j u(t)\|_p
\lesssim \int^t_0 \e^{-(c_02^{\alpha j}+\lambda-\lambda_0)(t-s)}\Big(2^{-\gamma j}
\big(\|u\|_{B^{1-\beta+\gamma}_{\infty,1}}+\|u\|_{B^{1-\beta+\gamma}_{p,1}}\big)
+\|\Pi_j f\|_p\Big)\dif s\no\\
&\quad\leq 2^{-\gamma j}\int^t_0 \e^{-(c_02^{\alpha j}+\lambda-\lambda_0)(t-s)}\Big(\|u\|_{B^{1-\beta+\gamma}_{\infty,1}}
+\|u\|_{B^{1-\beta+\gamma}_{p,1}}+\|f\|_{B^\gamma_{p,\infty}}\Big)\dif s\label{KeyE}\\
&\quad\leq 2^{-\gamma j}2^{\lambda_0 t}\int^t_0 \e^{-c_02^{\alpha j}s}\dif s
\Big(\|u\|_{L^\infty_t(B^{1-\beta+\gamma}_{\infty,1})}
+\|u\|_{L^\infty_t(B^{1-\beta+\gamma}_{p,1})}+\|f\|_{L^\infty_t(B^\gamma_{p,\infty})}\Big).\no
\end{align}
Hence,
\begin{align}\label{ED3}
&\|u(t)\|_{B^{\alpha+\gamma}_{p,\infty}}=\sup_{j\geq -1} \Big(2^{(\alpha+\gamma)j}\|\Pi_j u(t)\|_p\Big)
\lesssim\|u\|_{L^\infty_t(B^{1-\beta+\gamma}_{\infty,1})}
+\|u\|_{L^\infty_t(B^{1-\beta+\gamma}_{p,1})}+\|f\|_{L^\infty_t(B^{\gamma}_{p,\infty})},
\end{align}
where we have used that $2^{\alpha j}\int^t_0 \e^{-c_02^{\alpha j}s}\dif s=(1-\e^{-c_0 2^{\alpha j}})/c_0\leq 1/c_0$.
 
Let $\theta\in(0,\alpha+\beta -d/p-1)$. By embedding relation \eqref{Emm} and interpolation theorem, we have for all $\eps\in(0,1)$,
\begin{align*}
\|u\|_{B^{1-\beta+\gamma}_{\infty,1}}&\leq C\|u\|_{B^{\alpha+\gamma-\theta}_{p,\infty}}
\leq C\|u\|^{1-\frac{\theta}{\alpha}}_{B^{\alpha+\gamma}_{p,\infty}}\|u(s)\|^{\frac{\theta }{\alpha}}_{B^{\gamma}_{p,\infty}}
\leq \eps\|u\|_{B^{\alpha+\gamma}_{p,\infty}}+C_{\eps}\|u\|_{B^{\gamma}_{p,\infty}},
\end{align*}
and similarly,
$$
\|u\|_{B^{1-\beta+\gamma}_{p,1}}\leq \eps\|u\|_{B^{\alpha+\gamma}_{p,\infty}}+C_{\eps}\|u\|_{B^{\gamma}_{p,\infty}}.
$$
Substituting these into \eqref{ED3} and letting $\eps$ be small enough, we get
\begin{align}\label{EE4}
\|u\|_{L^\infty_t (B^{\alpha +\gamma}_{p,\infty})}\lesssim\|u\|_{L^\infty_t(B^{\gamma}_{p,\infty})}+\|f\|_{L^\infty_t(B^{\gamma}_{p,\infty})},
\end{align}
and also,
\begin{equation}\label{Midestimate}
\|u\|_{L^\infty_t(B^{1-\beta+\gamma}_{\infty,1})}+\|u\|_{L^\infty_t(B^{1-\beta+\gamma}_{p,1})}
\lesssim \|u\|_{L^\infty_t(B^{\gamma}_{p,\infty})}+\|f\|_{L^\infty_t(B^{\gamma}_{p,\infty})}.
\end{equation}
Now, multiplying both sides of \eqref{KeyE} by $2^{\gamma j}$ and then taking supremum over $j$, we obtain
\begin{align*}
\|u(t)\|_{B^{\gamma}_{p,\infty}}
&\lesssim\int^t_0 \e^{-(\lambda-\lambda_0)(t-s)}\Big(\|u\|_{B^{1-\beta+\gamma}_{\infty,1}}
+\|u\|_{B^{1-\beta+\gamma}_{p,1}}+\|f\|_{B^\gamma_{p,\infty}}\Big)\dif s\\
&\lesssim \int_0^t \|u\|_{L^\infty_s (B^\gamma_{p,\infty})} \dif s+\Big(1\wedge |\lambda-\lambda_0|^{-1}\Big)\|f\|_{L^\infty_t(B^{\gamma}_{p,\infty})}.
\end{align*}
Thus, by Gronwall's inequality we get
\begin{equation}\label{Loworder}
\|u\|_{L^\infty_T(B^{\gamma}_{p,\infty})}\leq \Big(1\wedge |\lambda-\lambda_0|^{-1}\Big)\|f\|_{L^\infty_T(B^{\gamma}_{p,\infty})},
\end{equation}
where $C=C(d,\a,\beta,p,T,\|b_1\|_{L^\infty_T(B^\beta_{p,\infty})},\|b_2\|_{L^\infty_T(B^\beta_{\infty,\infty})})$,
 which together with \eqref{EE4} yields \eqref{Hightorder}. 
Combining  \eqref{Hightorder} with \eqref{Loworder}, using the interpolation theorem again, we obtain \eqref{ED7}.

\medskip

(ii) Let $\rho$ be a non-negative smooth function with compact support in $\R^d$ and $\int_{\R^d} \rho(x)\d x=1$. 
Define $\rho_{\eps}(x):=\eps^{-d} \rho(\eps^{-1}x)$, $b_\eps:=\rho_\eps*b$, $f_\eps=:\rho_\eps*f$. 
Let $u_\eps$ be the smooth solution of PDE \eqref{EE8} corresponding to $b_\eps$ and $f_\eps$, which can be written as
\begin{align}\label{rep}
u_\eps(t,x)=\int^t_0\e^{-\lambda s}\mE\Big([b_\eps\cdot\nabla u_\eps+f_\eps](s,x+Z_{t-s})\Big)\dif s.
\end{align}
By the a priori estimate \eqref{Hightorder}, we have
$$
\sup_{0< \eps\leq 1}\|u_\eps\|_{L^\infty_T (B^{\alpha+\gamma}_{p,\infty})}\leq C \|f\|_{L^\infty_T(B^{\gamma}_{p,\infty})}.
$$
By this uniform estimate and \eqref{rep}, we also have for all $ 0\leq t\leq t'\leq T$, 
\begin{align*}
\lim_{|t-t'|\to0} \sup_{0<\eps\leq 1} \|u_\eps(t)-u_\eps(t')\|_\infty=0. 
\end{align*}
Thus, there is a subsequence (still denoted by $u_\eps$) and a continuous function $u$ with
$$
\|u\|_{L^\infty_T (B^{\alpha+\gamma}_{p,\infty})}\leq C \|f\|_{L^\infty_T(B^{\gamma}_{p,\infty})},~~\|u_\eps-u\|_{L^\infty([0,T];C^1(B_R))}\to 0.
$$
By taking limits and suitable weak convergence method, we obtain the existence of solution $u$ (see \cite{Ch-So-Zh} for more details).
\end{proof}

\subsection{Variable diffusion coefficient case}   

In this subsection we consider the variable diffusion coefficient case, 
and introduce the following assumptions on $\sigma_t (x)$:
\begin{enumerate}[{\bf (H$^\theta_\eps$)}]
\item There are $\theta,\eps\in(0,1)$ and $\Lambda\geq 1$ such that
\begin{align}
&\|\sigma_t(x)-\sigma_t(y)\|\leq \Lambda|x-y|^\theta \quad \hbox{and} \quad  
\sigma_t(x)=\sigma_t(0) \  \hbox{ for }    |x|\geq\eps,\label{DD1}\\
&\qquad\Lambda^{-1} |\xi|^2\leq|\sigma_t(0)\xi|^2\leq\Lambda|\xi|^2 
\quad \hbox{for every }  \xi\in\mR^d.      \label{DD3}
\end{align}
\end{enumerate}
Notice that \eqref{DD1} means that $\sigma$ only varies near $0$ and this implies that
\begin{align}\label{DD2}
\|\sigma_t(x+y)-\sigma_t(x)\|\leq C(|y|\wedge\eps)^\theta.
\end{align}
About the L\'evy measure $\nu$, we assume
\begin{enumerate}[{\bf (H$^\alpha_\nu$)}]
\item There are  $\nu_1,\nu_2\in\mL^{(\alpha)}_{non}$ so that 
$$
\nu_1(A)\leq\nu(A)\leq\nu_2(A) \quad \hbox{for every }  A\in\sB(B_1).
$$
In particular, for any $\delta>\alpha$, 
\begin{align}\label{Mom}
\int_{|z|\leq 1}|z|^\delta\nu(\dif z)<\infty.
\end{align}
\end{enumerate}
Since the L\'evy measure $\nu$ is not necessarily absolutely continuous with respect to the Lebesgue measure, it seems hard to show that
for any $f\in B^{\a+\gamma}_{p,\infty}$,
$$
x\mapsto \int_{|z|>1}f(x+\sigma_t(x) z)\nu(\dif z)\in B^\gamma_{p,\infty}, 
$$
which is very essential if  one wants to use the perturbation argument. 
Thus we have to first remove the large jump part and consider the following operator
$$
\widetilde{\sL_t}f(x):=\widetilde{\cL^\nu_{\sigma_t(x)}}f(x)=\int_{|z|\leq 1}\Big(f(x+\sigma_t(x) z)-f(x)-\sigma_t(x) z\cdot\nabla f(x)\Big)\nu(\dif z).
$$
The following theorem is the main result of this subsection. 
Although  this analytic result needs
a special assumption on the oscillation of $\si_t(\cdot)$, it is enough for us to get our Theorem \ref{Th11}.  

\bt\label{Th34}
Let $\beta\in(0,1)$ and $\alpha\in(0,2)$ with $\alpha+\beta>1$. Let $T>0$ and
$p\in(\frac{d}{\alpha+\beta-1}\vee \frac{d^2}{\alpha\wedge 1} \vee 2 ,\infty)$, 
$\theta\in( \beta,1]$.
Suppose {\bf (H$^{\alpha}_\nu$)} and $b=b_1+b_2$ with
$$
b_1\in L^\infty_T(B^\beta_{p,\infty}),\ \ b_2\in L^\infty_T(B^\beta_{\infty,\infty}).
$$ 
Then there are $\eps_0\in(0,1)$ and $\lambda_0>0$ such that for all $\eps\in(0,\eps_0)$ and
$\lambda\in(\lambda_0,\infty)$, under {\bf (H$^\theta_\eps$)}, for any  $\gamma\in(0,\beta]$ and $f\in L^\infty_T(B^\gamma_{p,\infty})$, 
there is  a unique solution $u\in L^\infty_T (B^{\alpha+\gamma}_{p,\infty})$ solving
$$
u(t)=\int^t_0[(\widetilde{\sL_s}-\lambda) u+b\cdot\nabla u+f]\dif s.
$$
Moreover, we have
\begin{align}\label{ED9}
\|u\|_{L^\infty_T (B^{\alpha+\gamma}_{p,\infty})}\leq C \|f\|_{L^\infty _T(B^{\gamma}_{p,\infty})},
\end{align}
and for any $s\in(0,\alpha+\gamma)$,
\begin{align}\label{ED91}
\|u\|_{L_T^\infty(B^s_{p,\infty})}\leq c_\lambda \|f\|_{L^\infty _T(B^{\gamma}_{p,\infty})},
\end{align}
where $c_\lambda=c(\lambda,d,p,\a,s)\to0$ as $\lambda\to \infty$. 
\et

In order to get the above result, we need the following commutator estimate.

\bl\label{Le34}
Under {\bf (H$^\theta_\eps$)} and {\bf (H$^\alpha_\nu$)}, 
for any $p>1$, we have
$$
\big\|[\Delta^{s/2},\widetilde{\sL_t}]u\big\|_p\leq 
C\left\{
\begin{aligned}
&\eps^{\theta\delta-s+d/p}\|u\|_{B^\delta_{\infty,\infty}},\ \, \alpha\in(0,1),\delta\in(\alpha,1],s\in(0,\theta\delta);\\
&\eps^{\theta-s+d/p}\|\nabla u\|_{B^{\delta-1}_{\infty,\infty}},\ \alpha\in[1,2),\delta\in(\alpha,2),s\in(0,\theta),
\end{aligned}
\right.
$$
where $[\Delta^{s/2},\widetilde{\sL_t}]u:=\Delta^{s/2}\widetilde{\sL_t}u-\widetilde{\sL_t}\Delta^{s/2}u$, and the constant $C>0$ is independent of $\eps$.
\el
\begin{proof}
We only prove it for $\alpha\in[1,2)$ since the case $\alpha\in(0,1)$ is similar. For simplicity of notation, we drop the variable $t$ in $\sigma_t(x)$ and write
\begin{align*}
\Gamma^\sigma_u(x,y,z)&:=u(x+y+\si(x+y)z)-u(x+y+\si(x)z)\\
&\qquad-(\sigma(x+y)-\sigma(x))z\cdot\nabla u(x+y),
\end{align*}
and
$$
\|[\Delta^{s/2},\widetilde{\sL_t}]u\|^p_p=\left(\int_{|x|\leq2\eps}+\int_{|x|>2\eps}\right)\big|[\Delta^{s/2},\widetilde{\sL_t}]u(x)\big|^p\dif x=:\cJ_1+\cJ_2.
$$
Let $\delta\in(\alpha,2)$. By \eqref{DD2} and \eqref{ED8}, we have
\begin{align*}
|\Gamma^\sigma_u(x,y,z)|&\lesssim(|y|\wedge\eps)^\theta |z|\int^1_0|\nabla u(x+y+(1-r)\si(x+y)z+r\si(x)z)\\
&\qquad\qquad\qquad\qquad-\nabla u(x+y)|\dif r\\
&\lesssim(|y|\wedge\eps)^\theta |z|^{\delta}\|\nabla u\|_{B^{\delta-1}_{\infty,\infty}},
\end{align*}
and by definition,
\begin{align*}
[\Delta^{s/2},\widetilde{\sL_t}]u(x)= \int_{|z|\leq 1}\nu(\d z) \int_{\R^d}\frac{\Gamma^\sigma_u(x,y,z)}{|y|^{d+s}}\d y.
\end{align*}
Thus, for $\cJ_1$, by \eqref{Mom} we have
\begin{align*}
\cJ_1&\lesssim\|\nabla u\|^p_{B^{\delta-1}_{\infty,\infty}}\int_{|x|\leq2\eps}
\left|\int_{|z|\leq 1}\nu(\d z) \int_{\R^d}\frac{|z|^{\delta}(|y|\wedge\eps)^\theta}{|y|^{d+s}}\d y\right|^p\dif x\\
&\lesssim\eps^d\|\nabla u\|^p_{B^{\delta-1}_{\infty,\infty}}
\left|\int_{|y|\leq\eps}\frac{|y|^\theta\dif y}{|y|^{d+s}}+\int_{|y|>\eps}\frac{\eps^\theta\dif y}{|y|^{d+s}}\right|^p
\lesssim\|\nabla u\|^p_{B^{\delta-1}_{\infty,\infty}}\eps^{(\theta-s)p+d}.
\end{align*}
For $\cJ_2$, since $\Gamma^\sigma_u(x,y,z)=0$ for $|x|, |x+y|>\eps$ by \eqref{DD1}, we have
\begin{align*}
\cJ_2&=\int_{|x|>2\eps}\left|\int_{|z|\leq 1}\nu(\d z) \int_{|x+y|\leq\eps}\frac{\Gamma^\sigma_u(x,y,z)}{|y|^{d+s}}\d y\right|^p\dif x\\
&\lesssim\|\nabla u\|^p_{B^{\delta-1}_{\infty,\infty}}\int_{|x|>2\eps}\left|\int_{|z|\leq 1}\nu(\d z) 
\int_{|x+y|\leq\eps}\frac{|z|^{\delta}(|y|\wedge\eps)^\theta}{|y|^{d+s}}\d y\right|^p\dif x\\
&\lesssim\|\nabla u\|^p_{B^{\delta-1}_{\infty,\infty}}\eps^{\theta p}\int_{|x|>2\eps}\left|
\int_{|x+y|\leq\eps}\frac{1}{|y|^{d+s}}\d y\right|^p\dif x\\
&\lesssim\|\nabla u\|^p_{B^{\delta-1}_{\infty,\infty}}\eps^{\theta p+dp}\int_{|x|>2\eps}\frac{1}{(|x|-\eps)^{p(d+s)}}\dif x
\lesssim\|\nabla u\|^p_{B^{\delta-1}_{\infty,\infty}}\eps^{(\theta-s) p+d}.
\end{align*}
Combining the above calculations, we obtain the desired estimate.
\end{proof}
\bl\label{Le35}
Under {\bf (H$^\theta_\eps$)} and {\bf (H$^\alpha_\nu$)}, for any
$p\in(\tfrac{d^2}{\alpha\wedge 1},\infty)$ and $\theta\in(\tfrac{d}{p(\alpha\wedge 1)},1)$,
we have
$$
\big\|(\widetilde{\cL^\nu_{\sigma_t(\cdot)}}-\widetilde{\cL^\nu_{\sigma_t(0)}})f\big\|_{B^{\gamma}_{p,\infty}}
\leq c_\eps
\left\{
\begin{aligned}
&\|f\|_{B^{\alpha+\gamma}_{p,\infty}},\ \alpha\in(0,1),\ \gamma\in(0,\tfrac{(p\alpha-d)\theta}{p(1-\theta)}\wedge \theta);\\
&\|f\|_{B^{\alpha+\gamma}_{p,\infty}},\ \alpha\in[1,2),\ \gamma\in(0,\theta),
\end{aligned}
\right.
$$
where $c_\eps\to 0$ as $\eps\to 0$.
\el
\begin{proof}
For simplicity of notation, we drop the time variable $t$ and write 
$$
\cT_\sigma f:=\widetilde{\cL^\nu_{\sigma(\cdot)}}f-\widetilde{\cL^\nu_{\sigma(0)}}f.
$$
We prove the estimate for $\alpha\in(0,1)$. The case $\alpha\in[1,2)$ is similar.
By \cite[(2.19)]{Ch-Zh2}, we have
\begin{align}
\|\widetilde{\cL^{\nu}_{\sigma_1}} f-\widetilde{\cL^{\nu}_{\sigma_2}} f\|_p
\lesssim(\|\sigma_1-\sigma_2\|^{\a}\wedge 1)\|f\|_{H^\alpha_p}.\label{ED5}
\end{align}
Noticing that by \eqref{DD1},
$$
|\cT_\sigma f(x)|\leq |\widetilde{\cL^\nu_{\sigma(x)}}f(x)-\widetilde{\cL^\nu_{\sigma(0)}}f(x)|\leq
\sup_{\|\sigma-\sigma(0)\|\leq\Lambda\eps^\theta}|\widetilde{\cL^\nu_{\sigma}}f(x)-\widetilde{\cL^\nu_{\sigma(0)}}f(x)|,
$$ 
since $p>d^2/\a$, by \cite[Lemma 2.2]{Ch-Zh2} and \eqref{ED5}, we have
\begin{align}\label{ED6}
\|\cT_\sigma f\|_p\leq \left\|\sup_{\|\sigma-\sigma(0)\|\leq\Lambda\eps^\theta}|\widetilde{\cL^\nu_{\sigma}}f(\cdot)-\widetilde{\cL^\nu_{\sigma(0)}}f(\cdot)|\right\|_p\lesssim
\eps^{\a\theta}\|f\|_{H^\alpha_p}.
\end{align}
To show the estimate, by \eqref{ED4}, we have
$$
\|\Pi_i\cT_\sigma f\|_p
\leq\sum_{j>i}\|\Pi_i\cT_\sigma\Pi_jf\|_p+\sum_{j\leq i}\|\Pi_i\cT_\sigma\Pi_jf\|_p=:\cJ_1+\cJ_2.
$$
For $\cJ_1$, by \eqref{ED6} and Bernstein's inequality \eqref{EE99}, we have
\begin{align*}
\cJ_1&\lesssim\sum_{j>i}\|\cT_\sigma\Pi_jf\|_p\lesssim \eps^{\a\theta}\sum_{j>i}\|\Pi_jf\|_{H^\alpha_p}
\lesssim \eps^{\a\theta}\sum_{j>i}2^{\alpha j}\|\Pi_jf\|_p\\
&\leq\eps^{\a\theta}\sum_{j>i}2^{-\gamma j}\|f\|_{B^{\alpha+\gamma}_{p,\infty}}
=\eps^{\a\theta}2^{-\gamma i}\|f\|_{B^{\alpha+\gamma}_{p,\infty}}/(1-2^{-\gamma}).
\end{align*}
For $\cJ_2$, for any $\gamma \in(0,\tfrac{(p\alpha-d)\theta}{p(1-\theta)}\wedge \theta)$, 
one can choose $\delta \in(\a,\frac{(p\alpha-d)}{p(1-\theta)}\wedge1]$ such that $\gamma<\delta \theta$.  
Since $\delta<\frac{(p\a-d)}{p(1-\theta)}$, we get $\delta+\frac{d}{p}-\a<\delta \theta$. By this,  
we can fix $s\in (\gamma, \delta\theta)$ such that $s>\delta+\frac{d}{p}-\a$. Using Bernstein's inequality and Lemma \ref{Le34}, we have 
\begin{align*}
\cJ_2&= \sum_{j\leq i}\|\Pi_i\Delta^{-s/2}\Delta^{s/2}\cT_\sigma\Pi_jf\|_p
\lesssim 2^{-s i} \sum_{j\leq i}\|\Delta^{s/2}\cT_\sigma\Pi_jf\|_p\\
&\leq 2^{-s i} \sum_{j\leq i}\left(\|[\Delta^{s/2},\cT_\sigma]\Pi_jf\|_p+\|\cT_\sigma\Delta^{s/2}\Pi_jf\|_p\right)\\
&= 2^{-s i} \sum_{j\leq i}\left(\|[\Delta^{s/2}, \widetilde{\sL_t} ]\Pi_jf\|_p+\|\cT_\sigma\Delta^{s/2}\Pi_jf\|_p\right)\\
&\lesssim 2^{-s i}\sum_{j\leq i}\left(\eps^{\theta\delta-s+d/p}\|\Pi_jf\|_{C^\delta}+\eps^{\a\theta} \|\Delta^{s/2}\Pi_jf\|_{H^\alpha_p}\right)\\
&\lesssim 2^{-s i}\sum_{j\leq i}\left(\eps^{\theta\delta-s+d/p}2^{(\delta+d/p)j}\|\Pi_jf\|_p+\eps^{\a\theta} 2^{(s+\alpha)j} \|\Pi_jf\|_p\right)\\
&\lesssim \eps^{\theta\delta-s+d/p}2^{-s i}\sum_{j\leq i} \left(2^{(\delta+d/p)j}+2^{(\alpha+s)j}\right)\|\Pi_jf\|_p\\
&\lesssim \eps^{\theta\delta-s+d/p} \left(2^{-si}\sum_{j\leq i}2^{(\delta+\frac{d}{p}-\a-\gamma)j} 
+2^{-s i}\sum_{j\leq i}2^{(s-\gamma)j}\right)\|f\|_{B^{\alpha+\gamma}_{p,\infty}}\\
&\lesssim \eps^{\theta\delta-s+d/p} 2^{-\gamma i}\|f\|_{B^{\alpha+\gamma}_{p,\infty}}.
\end{align*}
Combining the above calculations, we obtain the estimate.
\end{proof}

We are in a position to give
\begin{proof}[Proof of Theorem \ref{Th34}]
By the classical continuity method or Picard's iteration, it suffices to prove the a priori estimate \eqref{ED9}. Let us rewrite equation \eqref{EE8} as
$$
\partial_tu+\lambda u- \cL^\nu_{\sigma_t(0)}u-b\cdot\nabla u= 
f+(\widetilde{\cL^\nu_{\sigma_t(\cdot)}}-\widetilde{\cL^\nu_{\sigma_t(0)}})u+\int_{|z|>1}\Big(u(\cdot+\sigma_t(0) z)-u(\cdot)\Big)\nu(\dif z).
$$
By the assumption $\theta>\beta>1+\frac{d}{p}-\a$,  one sees that $\tfrac{(p\alpha-d)\theta}{p(1-\theta)}\geq \theta$. 
So, by \eqref{Hightorder} and Lemma \ref{Le35}, 
for any $\gamma\in(0,\beta]$, we have
\begin{align*}
c^{-1}_\lambda \|u\|_{L_T^\infty(B^\gamma_{p,\infty})}+\|u\|_{L_T^\infty(B^{\a+\gamma}_{p,\infty})}&\leq C \|f\|_{L^\infty_T(B^\gamma_{p,\infty})}
+\big\|(\widetilde{\cL^\nu_{\sigma_t(\cdot)}}-\widetilde{\cL^\nu_{\sigma_t(0)}})u\big\|_{L^\infty_T(B^\gamma_{p,\infty})}+C \|u\|_{L^\infty_T(B^{\gamma}_{p,\infty})}\\
&\leq C \|f\|_{L^\infty_T(B^\gamma_{p,\infty})}+ c_{\eps} \|u\|_{L^\infty_T(B^{\a+\gamma}_{p,\infty})}
+C \|u\|_{L^\infty_T(B^{\gamma}_{p,\infty})},
\end{align*}
where $C$ is independent of $\eps$ and $\lambda$.
Choosing $\eps_0$ small and $\lambda_0$ large enough, we get for all $\eps\in(0,\eps_0)$ and $\lambda\in(\lambda_0,\infty)$,
$$
c^{-1}_\lambda \|u\|_{L_T^\infty(B^\gamma_{p,\infty})}+\|u\|_{L_T^\infty(B^{\a+\gamma}_{p,\infty})}\leq C \|f\|_{L^\infty_T(B^\gamma_{p,\infty})}.
$$
Thus, we finish the proof.
\end{proof}

\section{Pathwise well-posedness of SDE \eqref{SDE}} 

In this section, we give a proof for the main result of this paper, Theorem \ref{Th11}.
Let $N(\dif t,\dif z)$ be the Poisson random measure associated with $Z$, that is,
$$
N((0,t]\times E)=\sum_{s\leq t}1_{E}(\Delta Z_s),\ E\in\sB(\mR^d),\ \Delta Z_s:=Z_s-Z_{s-},
$$
whose intensity measure is given by $\dif t\nu(\dif z)$. Let $\tilde N(\dif t,\dif z)=N(\dif t,\dif z)-\dif t\nu(\dif z)$ be the compensated Poisson random
martingale measure. By L\'evy-It\^o's decomposition, we have
$$
Z_t=\ell t+\int^t_0\!\!\!\int_{|z|\leq 1} z\tilde N(\dif s,\dif z)+\int^t_0\!\!\!\int_{|z|>1} zN(\dif s,\dif z),
$$
where $\ell :=\mE\left(Z_1-\int_{|z|>1}zN(1,\dif z)\right)$.
Thus, SDE \eqref{SDE} can be written as
\begin{align*}
X_t&=X_0+\int_0^t b^\ell_s(X_s)\d s+\int_0^t\!\!\!\int_{|z|\leq 1}\si_s(X_{s-})z\tilde N(\dif s,\dif z)+\int_0^t\!\!\!\int_{|z|>1}\si_s(X_{s-})zN(\dif s,\dif z), 
\end{align*}
where $b^\ell_s(x):=b_s(x)+\sigma_s(x)\ell$.
To solve SDE \eqref{SDE}, by standard interlacing technique, it suffices to solve the following SDE
\begin{equation}\label{SDE3}
X_t=X_0+\int_0^t b^\ell_s(X_s)\d s+\int_0^t\!\!\!\int_{|z|\leq 1}\si_s(X_{s-})z\tilde N(\dif s,\dif z).
\end{equation}
Below we shall fix a complete and right continuous filtered probability space $(\Omega,\sF,\mP; (\sF_t)_{t\geq 0})$ 
so that all the processes are defined on it. 

Let us first show the following preliminary result.

\bl\label{Le11}
Let $b_t(x,\omega)$ and $\si_t(x,\omega)$ be two $\sB(\mR_+)\times\sB(\mR^d)\times\sF_0$-measurable functions. Let $\beta\in(1-\tfrac{\alpha}{2},1)$ and 
$p\in(\frac{d}{\alpha/2+\beta-1}\vee \frac{d^2}{\alpha\wedge 1} \vee 2,\infty)$. Suppose that
$$
\sup_{\omega\in\Omega}\|b_\cdot(\cdot,\omega)\|_{L^\infty_T(B^\beta_{p,\infty})}<\infty,\ \ T>0,
$$
and $\sigma_\cdot(\cdot,\omega)$ satisfies {\bf (H$^1_{\eps_0}$)} with common bound $\Lambda$ for almost every $\omega$, where $\eps_0$ is the same as in Theorem \ref{Th34}.
Then there is a unique $\sF_t$-adapted solution $X_t$ so that
$$
X_t=X_0+\int_0^t b^\ell_s(X_s)\d s+\int_0^t\!\!\!\int_{|z|\leq 1} \si_s(X_{s-})z \tilde N(\dif s,\dif z).
$$ 
\el
\begin{proof}
Let $T>0$. 
Consider the following backward nonlocal parabolic system with random coefficients:
$$
\p_t\u_t+(\widetilde{\sL_t}-\lambda)\u_t+b^\ell\cdot\nabla \u_t+b=0,\ \ \u_T=0.
$$
By the assumptions and Theorem \ref{Th34}, for some $\eps_0\in(0,1)$ and $\lambda_0>0$, and for each $\omega$ and $\lambda>\lambda_0$,
there is a unique solution $\u_\cdot(\cdot,\omega)\in L^\infty_T(B^{\alpha+\beta}_{p,\infty})$ to the above equation with
$$
\sup_{\omega\in\Omega}\|\u_\cdot(\cdot,\omega)\|_{L_T^\infty(B^{\a+\beta}_{p,\infty})}\leq C,
$$
and for any $s\in(0,\alpha+\beta)$,
\begin{align}\label{ED99}
\sup_{\omega\in\Omega}\|\u_\cdot(\cdot,\omega)\|_{L_T^\infty(B^s_{p,\infty})}\leq c_\lambda,
\end{align}
where $c_\lambda\to 0$ as $\lambda\to\infty$.
In particular, one can choose $\lambda>\lambda_0$ large enough so that
\begin{align}\label{ED90}
\sup_{\omega\in\Omega}\|\nabla\u_\cdot(\cdot,\omega)\|_\infty\leq 1/2.
\end{align}
Since $\u$ is $\sF_0$-measurable, by It\^o's formula (cf. \cite{Ik-Wa}), we have 
\begin{align*}
\u_t(X_t)&=\u_0(X_0)+\int_0^t [\p_s \u_s+\widetilde\sL_s \u_s+b^\ell_s\cdot \nabla \u_s](X_s)\d s\\
&+\int_0^t\!\!\! \int_{|z|\leq 1} [\u_s(X_{s-}+\si_s(X_{s-})z)-\u_s(X_{s-})]\tilde{N} (\d s,\d z).
\end{align*}
Let $\Phi_t(x,\omega)=x+\u_t(x,\omega)$. Then by \eqref{ED90},
$x\mapsto\Phi_t(x,\omega)$ is a $C^1$-diffeomorphism and
\begin{equation}
\begin{split}
\label{Y1}
Y_t&:=\Phi_t(X_t)=\Phi_0(X_0)+\int_0^t\Big( \lambda \u_s(X_s)+\sigma_s(X_s)\ell\Big)\d s\\
&+\int_0^t\!\!\!\int_{|z|\leq 1} [\Phi_s(X_{s-}+\si_s(X_{s-})z)-\Phi_s(X_{s-})]\tilde{N} (\d s,\d z)\\
&=\Phi_0(X_0)+\int_0^t a_s(Y_s)\d s+\int_0^t\!\!\!\int_{|z|\leq 1} g_s(Y_{s-},z)\tilde{N} (\d s,\d z),
\end{split}
\end{equation}
where 
$$
a_t(y):=(\lambda \u_t+\sigma_t\ell)(\Phi^{-1}_t(y)),\ \ g_t(y,z):=\Phi_t(\Phi_t^{-1}(y)+\si_t(\Phi_t^{-1}(y))z)-y.
$$ 
Fix $\eta\in(\alpha/2,\alpha+\beta-1-d/p)$. Noticing that
$$
|[f(x+z)-f(x)]-[f(y+z)-f(y)]|\leq \|\nabla f\|_{B^\eta_{\infty,\infty}} |x-y|\,|z|^\eta, 
$$ 
we have
\begin{equation*}
\begin{split}
|g_t(x,z)-g_t(y, z)|
&=\big|[\u_t(\Phi_t^{-1}(x)+\si_t(\Phi^{-1}_t(x))z)-\u_t(\Phi_t^{-1}(x))+\si_t(\Phi^{-1}_t(x)z]\\
&-[\u_t(\Phi_t^{-1}(y)+\si_t(\Phi^{-1}_t(y))z)-\u_t(\Phi_t^{-1}(y))+\si_t(\Phi^{-1}_t(y)z]\big|\\
&\leq \big|[\u_t(\Phi_t^{-1}(x)+\si_t(\Phi^{-1}_t(x))z)-\u_t(\Phi_t^{-1}(x))]\\
&\qquad-[\u_t(\Phi_t^{-1}(y)+\si_t(\Phi^{-1}_t(x))z)-\u_t(\Phi_t^{-1}(y))]\big|\\
&+\big|\u_t(\Phi_t^{-1}(y)+\si_t(\Phi^{-1}_t(x))z)-\u_t(\Phi_t^{-1}(y)+\si_t(\Phi^{-1}_t(y))z)\big|\\
&\qquad+\big|\si_t(\Phi^{-1}_t(x)z-\si_t(\Phi^{-1}_t(y)z\big|\\
&\lesssim\big|\Phi^{-1}_t(x)-\Phi^{-1}_t(y)\big|\Big(\|\nabla\u_t\|_{B^\eta_{\infty,\infty}}|z|^\eta+\|\nabla \u_t\|_\infty|z|+\Lambda|z|\Big).
\end{split}
\end{equation*}
Thus, by \eqref{ED99} and Sobolev's embedding, we get for all $x,y\in\mR^d$ and $|z|\leq 1$,
$$
|g_t(x,z)-g_t(y, z)|\lesssim |x-y|\,|z|^\eta . 
$$
Moreover, we also have
$$
|a_t(x)-a_t(y)|\lesssim |x-y|.
$$
Since the coefficients of SDE \eqref{Y1} are Lipschitz continuous, 
by the  classical result, SDE \eqref{Y1} admits a unique solution (cf. \cite{Ik-Wa}). In particular, one can check that
$X_t=\Phi^{-1}_t(Y_t)$ satisfies the original equation \eqref{SDE3}. 
The proof is complete.
\end{proof}

We also need the following technical lemma in order to patch up the solution.

\bl\label{Le42}
Let $X_t$ be a $\mR^d$-valued right continuous process. Let $\tau$ be an $\sF_t$-stopping time. 
Suppose that for each $t\geq 0$,
$X_{t+\tau}$ is $\sF_{t+\tau}$-measurable. Then for each $t\geq 0$, $1_{\{\tau\leq t\}}X_t$ is $\sF_t$-measurable.
\el

\begin{proof}
Since $X_t$ is right continuous, we have
\begin{align*}
1_{\{\tau\leq t\}}X_t=\lim_{n\to\infty}1_{\{\tau\leq t\}}X_{t+\tau-[2^n\tau]2^{-n}}
=\lim_{n\to\infty}\sum_{j=0}^{[2^nt]}1_{\{\tau\leq t\}}  X_{t+\tau-j2^{-n}} 1_{\{j\leq 2^n\tau<j+1\}}.
\end{align*}
On the other hand, since by assumption $X_{t+\tau-j2^{-n}} = X_{\tau + (t-j2^{-n})} $ is $\sF_{\tau+(t -j2^{-n})}$-measurable
and $\tau+t-j2^{-n}$ is a stopping time 
when $t\geq  j2^{-n} $, we have for each $n\geq 1$ and $1\leq j\leq [2^nt]$, 
\begin{align*}
1_{\{\tau\leq t\}}X_{t+\tau-j2^{-n}}1_{\{j\leq 2^n\tau<j+1\}} 
 &= X_{t+\tau-j2^{-n}}  1_{\{2^{-n}j\leq \tau< (j+1) 2^{-n}\}} 1_{\{\tau\leq t\}}\\
&=X_{\tau+t-j2^{-n}}  1_{\{\tau+t-j2^{-n}<t+2^{-n}\}} 1_{\{j2^{-n}\leq\tau\leq t\}}. 
\end{align*}
 Noticing $X_{\tau+t-j2^{-n}}  1_{\{\tau+t-j2^{-n}<t+2^{-n}\}} \in \sF_{t+2^{-n}}$ and $1_{\{j2^{-n}\leq \tau\leq t\}}\in \sF_{t}$, we get
$$1_{\{\tau\leq t\}}X_t\in\cap_{n\geq 1 }\sF_{t+2^{-n}}=\sF_t.$$ 
The proof is complete.
\end{proof}
Now we can give

\begin{proof}[Proof of Theorem \ref{Th11}]
 By the discussion at the beginning of this section, we only need to prove the global well-posedness of \eqref{SDE3}. 
 Noticing $\E|X_t|^2\leq |X_0|^2+\|b\|^2_\infty t^2+\|\sigma\|_\infty^2 \int_{|z|\leq 1}  |z|^2\nu(\d z) t\leq C(1+t^2)$,  we get $\zeta=\infty$.  By Remark \ref{Rek1.3}, we can further assume that $b$ has support contained in ball $B_R$. 
Let $p\geq 1$. By definition \eqref{EE2}, we have
\begin{align*}
\|\Pi_j b\|_p^p=&\int_{\R^d} \left|\int_{B_R} h_j(x-y) b(y) \dif y\right|^p \d x
\leq    \|\Pi_j b\|_\infty^p|B_{2R}|+ \|b\|^p_\infty\int_{B^c_{2R}} \left(\int_{B_R} |h_j(x-y)|\dif y\right)^p \d x.
\end{align*}
 Noting that $h_j(x)=2^{jd}h_0(2^{j}x)$ by \eqref{EE7}, we have 
\begin{align*}
&\int_{B^c_{2R}} \left(\int_{B_R} |h_j(x-y)|\dif y\right)^p \d x
\leq \|h_0\|_1^{p-1} \int_{B^c_{2R}}\!\int_{B_R} |h_j(x-y)|\dif y\d x \\
 &\qquad\leq C({h_0}) \int_{B^c_{2R}}\!\int_{B_R}2^{jd}(2^j|x-y|)^{-2d}\dif y\d x \leq C({h_0,d, R}) 2^{-jd},
\end{align*}
where the second inequality is due to the polynomial decay property of $h_0$.
Hence,  
$$
\|b\|_{B^\beta_{p,\infty}}=\sup_{j\geq -1} 2^{j\beta} \|\Pi_j b\|_p \leq C \sup_{j\geq -1}2^{j\beta}(\|\Pi_j b\|_\infty+ 2^{-jd} \|b\|_\infty)\leq C\|b\|_{C^\beta}. 
$$
 Let $X_t$ be a solution of SDE \eqref{SDE}. Fix $\eps\in(0,1)$ being small as in Lemma \ref{Le11}. 
Let $\tau_0=0$. We define a sequence of stopping times as follows: for $n\in\mN$,
$$
\tau_n:=\inf\{t>\tau_{n-1}: |X_t-X_{\tau_{n-1}}|\geq\eps/2\}.
$$
We use induction to show the strong well-posedness of SDE \eqref{SDE}.
Suppose that  we have shown the existence and uniqueness of solutions up to time $\tau_n$. That is, there is a unique solution $X_t$
satisfying
$$
X_t=X_0+\int^t_0b^\ell_s(X_s)\dif s+\int^t_0\!\!\!\int_{|z|\leq 1}\sigma_s(X_{s-})z\tilde N(\dif s,\dif z),\ \ t<\tau_n,
$$
where $b^\ell_s(x):=b_s(x)+\sigma_s(x)\ell$.
Now, define
$$
X'_0:=X_{\tau_n}:=X_{\tau_n-}+\sigma_{\tau_n}(X_{\tau_n-})\Delta Z_{\tau_n}
$$
and
$$
\sF'_t:=\sF_{t+\tau_n},\ \ t\geq 0.
$$
Clearly, $X'_0\in\sF'_0$. We also introduce $\sF'_0$-measurable function $b'$ and $\sigma'$ as follows:
$$
b'_t(x,\omega)=b^\ell_{t+\tau_n}(x+X'_0(\omega))
$$
and 
$$
\si'_t(x,\omega):= 
\left\{
\begin{array}{ll}
{\si}_{t+\tau_n}(x+X'_0(\omega)), & |x|\leq \eps/2, \\
\left(
\begin{aligned}
&\tfrac{2(\eps-|x|)}{\eps}{\si}_{t+\tau_n}\big(\tfrac{\eps x}{2|x|}+X'_0(\omega)\big)\\
&+\tfrac{2(|x|-\eps/2)}{\eps}\si_{t+\tau_n}\big(X'_0(\omega)\big)
\end{aligned}
\right),&\eps/2<|x|\leq \eps,\\
\si_{t+\tau_n}(X'_0(\omega)), & |x|>\eps. 
\end{array}
\right.
$$
It is easy to see that $\si'$ satisfies {\bf (H$^1_{\eps}$)}.
Thus, by Lemma \ref{Le11}, the following SDE admits a unique solution
$$
X'_t=\int^t_0b'_s(X'_s)\dif s+\int^t_0\!\!\!\int_{|z|\leq 1}\sigma'_s(X'_{s-})z\tilde N(\tau_n+\dif s,\dif z).
$$
Define $\sF'_t$-stopping time
$$
\tau':=\inf\{t>0: |X'_t|\geq\eps/2\},\ \tau_{n+1}=\tau'+\tau_n,
$$
and for $t\in[0,\tau')$,
$$
X_{t+\tau_n}(\omega):=X'_{t}(\omega)+X_{\tau_n}(\omega).
$$
Observe that
$$
\{\tau_{n+1}<s\}=\cup_{t\in\mQ, t<s}(\{\tau'<t\}\cap\{\tau_n<s-t\})\in\sF_s.
$$
This means that $\tau_{n+1}$ is an $\sF_t$-stopping time. Since $t\mapsto X'_t$ is right continuous, by Lemma \ref{Le42}, we also have
$$
1_{[\tau_n,\tau_{n+1})}(t)X_t\mbox{ is $\sF_t$-measurable.}
$$
Moreover, by construction and induction, it is easy to see that for $t<\tau_{n+1}$, $X_t$ uniquely solves
$$
X_t=x+\int^t_0b_s(X_s)\dif s+\int^t_0\!\!\!\int_{|z|\leq 1}\sigma_s(X_{s-})z\tilde N(\dif s,\dif z).
$$
In fact, for $t\in[0,\tau_{n+1}-\tau_n)$, we have
\begin{align*}
X_{t+\tau_n}&=X_{\tau_n}+\int^t_0b_{s+\tau_n}(X_{s+\tau_n})\dif s+\int^{t}_{0}\!\!\!\int_{|z|\leq 1}\sigma_{s+\tau_n}(X_{(s+\tau_n)-})z\tilde N(\tau_n+\dif s,\dif z)\\
&=X_{\tau_n}+\int^{t+\tau_n}_{\tau_n}b_s(X_s)\dif s+\int^{t+\tau_n}_{\tau_n}\!\!\!\int_{|z|\leq 1}\sigma_s(X_{s-})z\tilde N(\dif s,\dif z)\\
&=x+\int^{t+\tau_n}_0b_s(X_s)\dif s+\int^{t+\tau_n}_0\!\!\!\int_{|z|\leq 1}\sigma_s(X_{s-})z\tilde N(\dif s,\dif z).
\end{align*}
The proof is complete.
\end{proof}

\end{document}